\documentclass[conference]{IEEEtran}
\IEEEoverridecommandlockouts
\usepackage{cite}
\usepackage{amsmath,amssymb,amsfonts}
\usepackage{algorithm}
\usepackage{algorithmic}
\usepackage{graphicx}
\usepackage{textcomp}
\usepackage{graphics}
\usepackage{subfig}
\usepackage{color, colortbl}
\usepackage[utf8x]{inputenc}
\def\BibTeX{{\rm B\kern-.05em{\sc i\kern-.025em b}\kern-.08em
    T\kern-.1667em\lower.7ex\hbox{E}\kern-.125emX}}
\begin{document}

\title{Optimized tour planning for drone-based urban traffic monitoring}

\author{C. Christodoulou, and P. Kolios
	\thanks{The authors are with the KIOS Research and Innovation Centre of Excellence (KIOS CoE) and the Department of Electrical and Computer Engineering, University of Cyprus, Nicosia, 1678, Cyprus. {\tt\small \{christodoulou.chrystalleni, pkolios\}@ucy.ac.cy}}
}
\maketitle

\begin{abstract}
Drones or Unmanned Aerial Vehicles (UAVs) have become a reliable and efficient tool for road traffic monitoring. Compared to loop detectors and bluetooth receivers (with high capital and operational expenditure), drones are a low-cost alternative that offers great flexibility and high quality data. 

In this work, we derive optimized tour plans that a fleet of drones can follow for rapid traffic monitoring across particular regions of the transportation network. To derive these tours, we first identify monitoring locations over which drones should fly through and then compute minimum travel-time tours based on realistic resource constrains. Evaluation results are presented over a real road network topology to demonstrate the applicability of the proposed approach.
\end{abstract}


\section{INTRODUCTION}\label{sec:intro}
UAVs are currently being used in a wide range of applications concerning large infrastructures and the environment \cite{b1}. This rapidly improving technology can be especially beneficial for fast and reliable information gathering where the deployment of fixed sensors is not reasonably possible \cite{DronesAugmentingOurLife}.  

In the case of road transportation networks, the advantage of drones is that they are not restricted to travelling over the road network and thus can swiftly move over disperse locations to capture road traffic data. In addition, compared to traditional sensors for traffic monitoring, UAVs can provide significantly better information over tighter time periods \cite{b2}. As a result they can provide valuable information more efficiently and faster. 

Previous research has concentrated mostly on the challenges of remotely acquiring data, using a single UAV with usually a fixed trajectory with the focus being mainly on deriving the process of identifying and tracking vehicles (\cite{b3} - \cite{b9}). Solutions with multiple UAVs and adaptive trajectories have also been studied with the aim to also improve detection and tracking performance (\cite{b10}, \cite{b11}). 

In this work we are concerned with the problem of identifying monitoring locations and deriving optimized tours to visit those locations in the least amount of time which is a necessary step prior to detecting and tracking vehicles. Our previous work in \cite{b12} discussed the image-processing related aspects that govern the optimal selection of monitoring locations while in this paper we extend that architecture to derive optimized tour plans.

The main parameters affecting monitoring locations include the line-of-sight obstacles between the roadways and the drone locations (including the building heights as well as the road and pavement widths) and the onboard camera field of view as a function of the drone altitude. Under this setting, the aim is to compute monitoring locations and derive the minimum-cost (in terms of distance and time) tours to be traversed by drone units in order to cover a certain area of interest. 

The computed tours are constraint by the maximum allowed flight distance (related to the limited fly time) of the fleet of drones in use. The resulting problem is formulated as a variant of the Multiple Travelling Salesman Problem (mTSP) which is an NP-Hard problem. Hence an algorithmic approach is developed based on k-means clustering and the Cheapest Insertion Algorithm (CIA) which is a well known heuristic for the TSP \cite{b13}. 

The rest of the paper is organized as follows. Related  work  is  reviewed in  Section \ref{related} where our proposed novelties and contributions are discussed. Section \ref{model} describes the deployment strategy that determines drone locations for effective traffic monitoring. In Section \ref{design} the tour planning algorithm is presented. Section \ref{sec:evaluation} provides evaluation results for a case study conducted with real data in the capital city of Cyprus, Nicosia. Finally, the paper is concluded in Section \ref{concl} with the key contributions and findings of the proposed work.   

\section{RELATED WORK}\label{related}
To obtain sufficient road traffic information, measurement data should be collected over an adequately large area and an adequately long period of time in order to be able to observe the dynamic behaviour of the underlying road network. Traditional technology tools for traffic monitoring, including loop detectors or static video cameras, are positioned at very specific locations in the transportation network. Hence the collected measurements are inherently very localized and fail to capture realistic mobility patterns including accurate vehicle trajectories, driving patterns, and individualized routes through the network. Drones on the other hand have the advantage of being both mobile, and able to capture data over an extended space-time dimension. 

The works in \cite{b14} and \cite{b15} analyzed the collection of video data using camera-equipped drones for traffic management. Capitalizing on the fact that aerial video feeds cover large stretches of road, the authors discuss the potential of using this information for a plethora of transportation operations, including: emergency vehicle guidance, track vehicle movements, estimation of typical roadway and parking area usage. 

Several additional works have looked into the problem of traffic monitoring using drone-based systems \cite{b3}-\cite{b11}. Significant research work is also placed on the video processing aspects on monitoring traffic by designing highly accurate and reliable vehicle detectors, tracking algorithms and processes to extract traffic parameters from the derived trajectories \cite{b2},\cite{b16}.

Our previous work \cite{b12} extends these approaches by incorporating data collection and processing capabilities from drone footage into a holistic framework that considers the drone capacities as well as the limitations presented by the underlying road traffic network to implement optimization strategies for drone-based traffic monitoring. As emphasized in the introduction, our contributions focus on deriving tour plans that drone units should pass through, in order to capture traffic patterns above an arbitrary stretch of a transportation network. 

\section{SYSTEM MODEL}\label{model}
As introduced in Section \ref{sec:intro}, our aim is to monitor a specific area of the road network by a fleet of drones. To do that, we first discretize the area into a set of locations out of which a subset of points will be used as monitoring locations and through which optimized tour plans will be computed. 

The set of locations (hereafter termed Points of Interest - PoI) is constructed by exploiting the fact that individual subareas of the region of interest are limited to a maximum building height $h_{M}$, above which no infrastructure in the particular subarea can exceed. Then there is a certain height $h_{U}$ that ensures line-of-sight for all vehicles located in a circular region of radius $R$ around those positions above the ground. Having $R_l$ for each subarea $l$, enables us to construct a non-uniform grid of locations by considering points on a square lattice with horizontal/vertical distance $\frac{R_l}{\rho}$, where $\rho$ is a granularity constant. Hence, for denser sampling the value of $\rho$ should be higher; resulting to a denser grid and thus a larger number of monitoring locations.
\begin{figure}
	\centering
	\includegraphics[width=\linewidth]{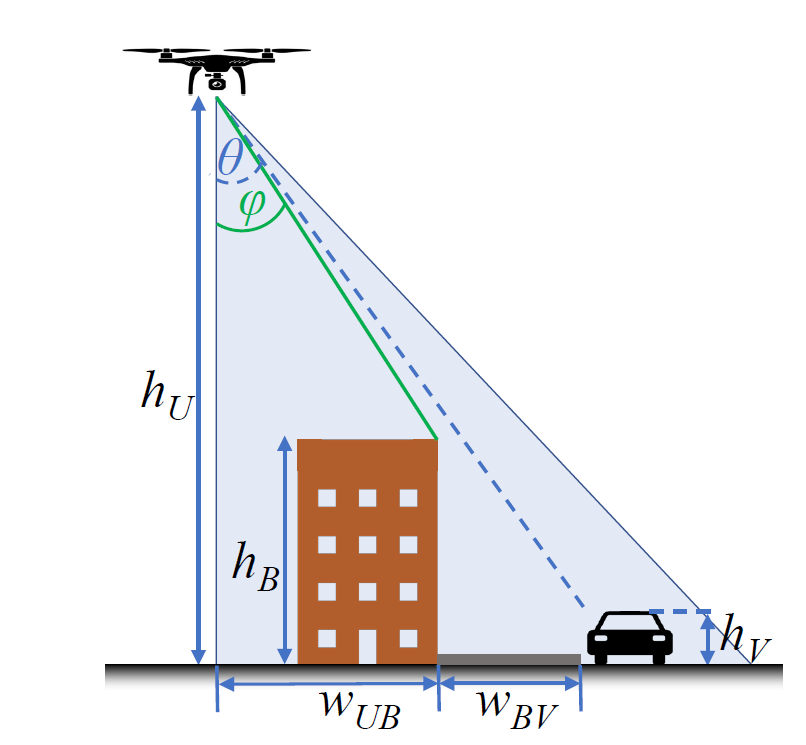}
	\caption{Geometric representation of the line-of-sight between a drone and a vehicle inside a transportation network.}
	\label{fig:UAV}
\end{figure}

To aid understanding in calculating height $h_{U}$ and radius $R$ for each subarea, lets first consider Fig. \ref{fig:UAV} showing a single drone over some arbitrary road segment. Let $w_{UB}$ and $w_{BV}$ denote the horizontal distance between the drone and the building that possibly obstructs the line-of-sight with passing vehicles, and the distance between that building and the vehicles, respectively, as shown in the figure. Consider the height of an arbitrary vehicle $h_{V}$ passing through the road segment, then to be able to observe the vehicle from the overlaying camera, a line-of-side should be achieved between the drone and the vehicle. As shown in the figure, this is only possible when $\theta \geq \phi$, where
\begin{align} 
	tan{(\theta)}&=\frac{w_{UB}+w_{BV}}{h_{U}-h_{V}}\\
	tan{(\phi)}&=\frac{w_{BV}}{h_{B}-h_{V}}
\end{align} 

The maximum radius at which line-of-sight is maintained is when $\theta = \phi$, hence $R$ can be computed as follows:
\begin{equation}
R=\frac{h_{U}-h_{V}}{h_{M}-h_{V}} w_{BV}
\end{equation}

Let $N$ and $M$ denote the sets of all potential drone locations and PoIs of the road network, respectively. Then, a drone positioned at $i \in N$ can monitor vehicles located along PoI $j \in M$ if $d_{ij} \leq R_{j}$, where $d_{ij} = \sqrt{(x_{i} - x_{j})^2 + (y_{i} - y_{j})^2}$, $R_{j}$ is the maximum detection range of PoI $j$ from a drone, while ($x_{i} , y_{i}$) and ($x_{j} , y_{j}$) are the coordinates of points $i$ and $j$, respectively.

Let matrix $C \in {\{0,1\}^{|N|\times|M|}}$contain entries $C_{ij}$ denoting whether PoI $j \in M$ can be monitored from a drone located at point $i\in N$ (i.e., $C_{ij} = 1$) or not (i.e., $C_{ij} = 0$). Then  decision variable $x_{i}$, $i \in N$ indicates whether a drone is located at $i$ (i.e., $x_{i} = 1$) or not (i.e., $x_{i} = 0$). Then, the placement problem can be solved using formulation $(P1)$ depicted below:
\begin{align}
(P1)& \>min \sum_{i \in N} x_{i}&\\
&s.t. \>\sum_{i \in N} C_{ij}x_{i} \geq 1, & j \in M \\
& x_{i} \in \{0,1\}, & i \in N
\end{align}

Problem (P1) aims at minimizing the number of drone locations $N$ needed to cover all PoIs. The solution to this problem provides a map with assigned locations that need to be visited to monitor the complete road network. As a note, $(P1)$ is a set covering problem which is known to be NP-hard \cite{b18}. However, this problem does not need to be solved in an online fashion and thus standard Mixed Integer Linear Program (MILP) solvers can be employed to solve it (including Gurobi \cite{b19} that we employ in Section \ref{sec:evaluation}).

\section{Computing tour plans}\label{design}
Given the set of locations to be visited, the aim here is to compute optimized tour plans that minimize the number of drones necessary to cover the target area, considering the limited fly times of the fleet of drones. To simplify the problem we assume that all drones have similar characteristics, fly at constant speed $v$, and with the same maximum fly time $t$ and the maximum travel distance is $m_d$.

\subsection{Cheapest Insertion Algorithm (CIA)}
\label{sec:CIA}
Let $G=(V,E)$ be a complete edge-weighted undirected graph, with $V=\{1,...,N\}$. Let the cost $c_{ij}$ indicate the length of an arbitrary edge $(i,j)\in E$. Then, in the basic Cheapest Insertion Algorithm, an initial tour $T$ is created over graph $G$ by any two nodes $i_{1}$ and $i_{2}$ such that:
\begin{equation}
c_{i_{1}}+c_{i_{2}}= \min_{i,j\in V, i\neq j} {c_{ij}+c_{ji}}
\end{equation}   

That is, the shortest two-node subtour is chosen as the initial tour $T$ with cost $\mathcal{C}_{0}=c_{i_{1}}+c_{i_{2}}$ . Thereafter, in an iterative manner a node $k$ is inserted into $T$ that has the cheapest cost $c_{k}$ from all nodes not in the subtour $T$. When a new node $k$ is added in the subtour between nodes $i$ and $j$, the length of the subtour increases by $c_{ik}+c_{kj}-c_{ij}$ while the new cost for every other node $k$ not in the subtour is calculated by:
\begin{equation}
c_{k}= \min_{(i,j)\in T} {c_{ik}+c_{kj}-c_{ij}}
\end{equation} 

Hence, the cost of the new subtour in each iteration is calculated as follows:
\begin{equation}
\mathcal{C}_{n}= \mathcal{C}_{n-1} + \min_{k\notin T} {c_{k}}
\end{equation} 

CIA is a greedy heuristic which at each iteration adds the node that increases the length of the current subtour as little as possible, to compute a minimum cost tour that passes through all nodes with the minimum cost. In our case we want to create multiple such tours with a finite total tour length.
\begin{figure}
	\centering
	\includegraphics[width=3.2in]{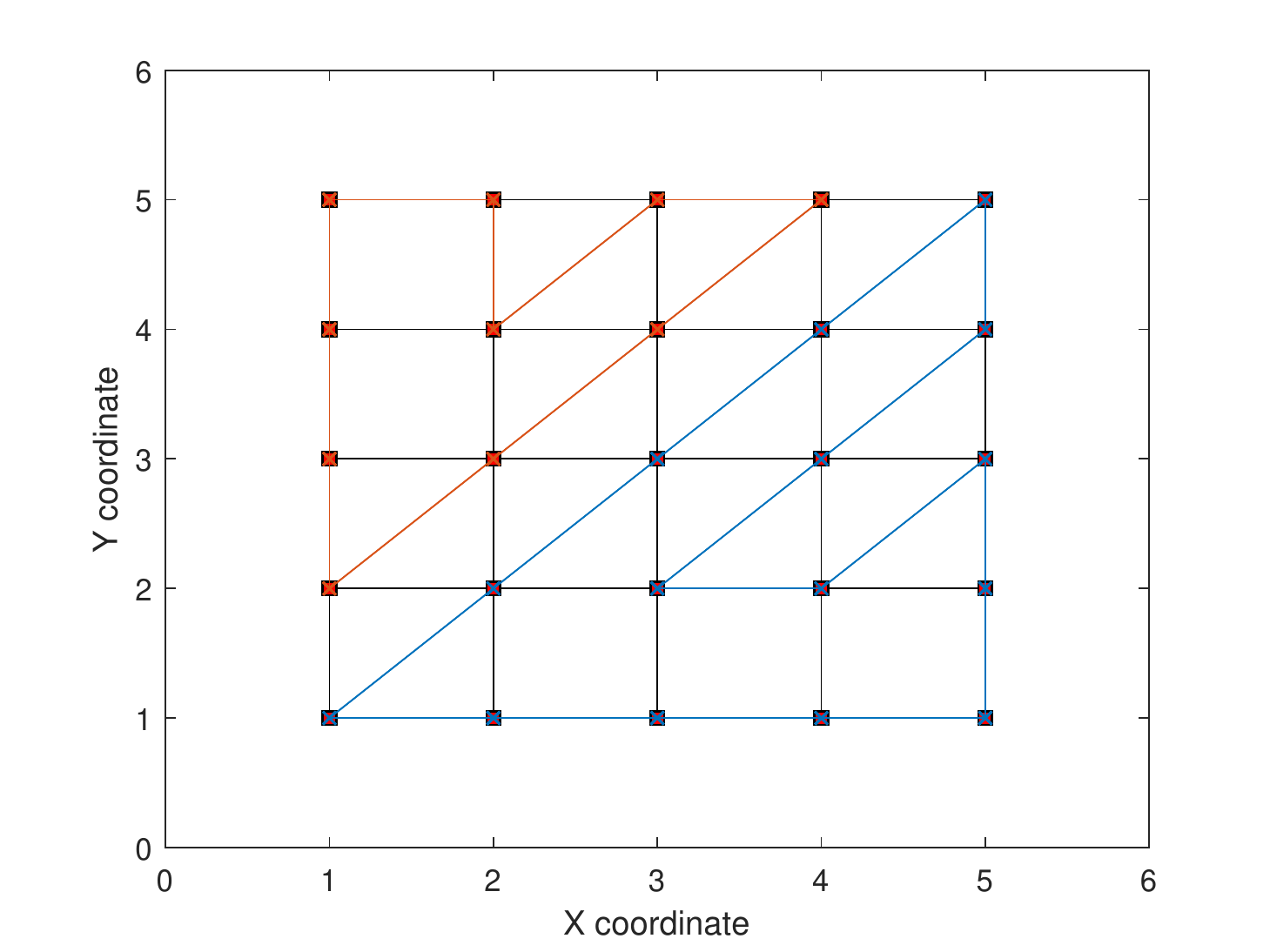}
	\caption{Example paths created using CIA.}
	\label{fig:CIA}
\end{figure}

Fig. \ref{fig:CIA} illustrates two such paths (with red and blue) created an a $5\times 5$ grid graph using CIA for two arbitrary subsets of nodes. Evidently, in order for all nodes to be visited under finite tour costs then an appropriate number of clusters must be found for which the tours created are feasible. 

\subsection{Multiple Tour Algorithm (MTA)}
As indicated above, the problem that arises is to compute the minimum number of tours that pass through the monitoring locations under some finite tour cost (i.e., the maximum travel distance that each drone can cover or a maximum revisit time).

In this work, we address this problem by developing an iterative two-stage algorithm as follows. The first stage splits the set of monitoring locations in $\mathcal{L}=\{1,\dots,L\}$ clusters. Then for each cluster CIA is employed (as described above) to compute a tour. The process repeats until all computed tours are feasible (i.e. the total travel cost of the computed tour ($\mathcal{C}_{l}$) is below the maximum acceptable tour length).

Since the main parameter affecting the tour length is the travel distance then the clustering is done by computing the distance between each monitoring location to each cluster centre, and then classifying the point to be in the cluster whose centre is closest to it. In this work, we employ k-means clustering to compute the node subsets.
\begin{algorithm}
	\begin{algorithmic}[1]
	\REQUIRE $L=1$
	\STATE Solve $(P1)$
	\WHILE{$\exists l\in\mathcal{L}$ such that $\mathcal{C}_{l}>m_d$}
	\STATE Create $L$ clusters using K-means
	\STATE Employ CIA described in Section \ref{sec:CIA} to develop a tour for each cluster
	\STATE $L=L+1$
	\ENDWHILE
	\end{algorithmic}
\caption{Multiple Tour Algorithm}
\label{alg:1}
\end{algorithm}

As detailed in Alg. \ref{alg:1}, steps 2-6 repeat until all computed tours that visit the graph nodes are feasible with respect to the total travel cost. Otherwise, the number of clusters required incrementally increases and the process repeats.

\section{EVALUATION RESULTS}\label{sec:evaluation}
To evaluate the performance of the proposed algorithm, we consider the metropolitan area of Nicosia in Cyprus as our test scenario. The road network, extracted from OpenStreetMap (www.openstreetmap.org), is comprised of all primary, secondary and residential links, and spans an area of approximately 100 $km^2$. 
\begin{figure*}[!t]
	\centering
	\subfloat[]{\includegraphics[width=2.4in]{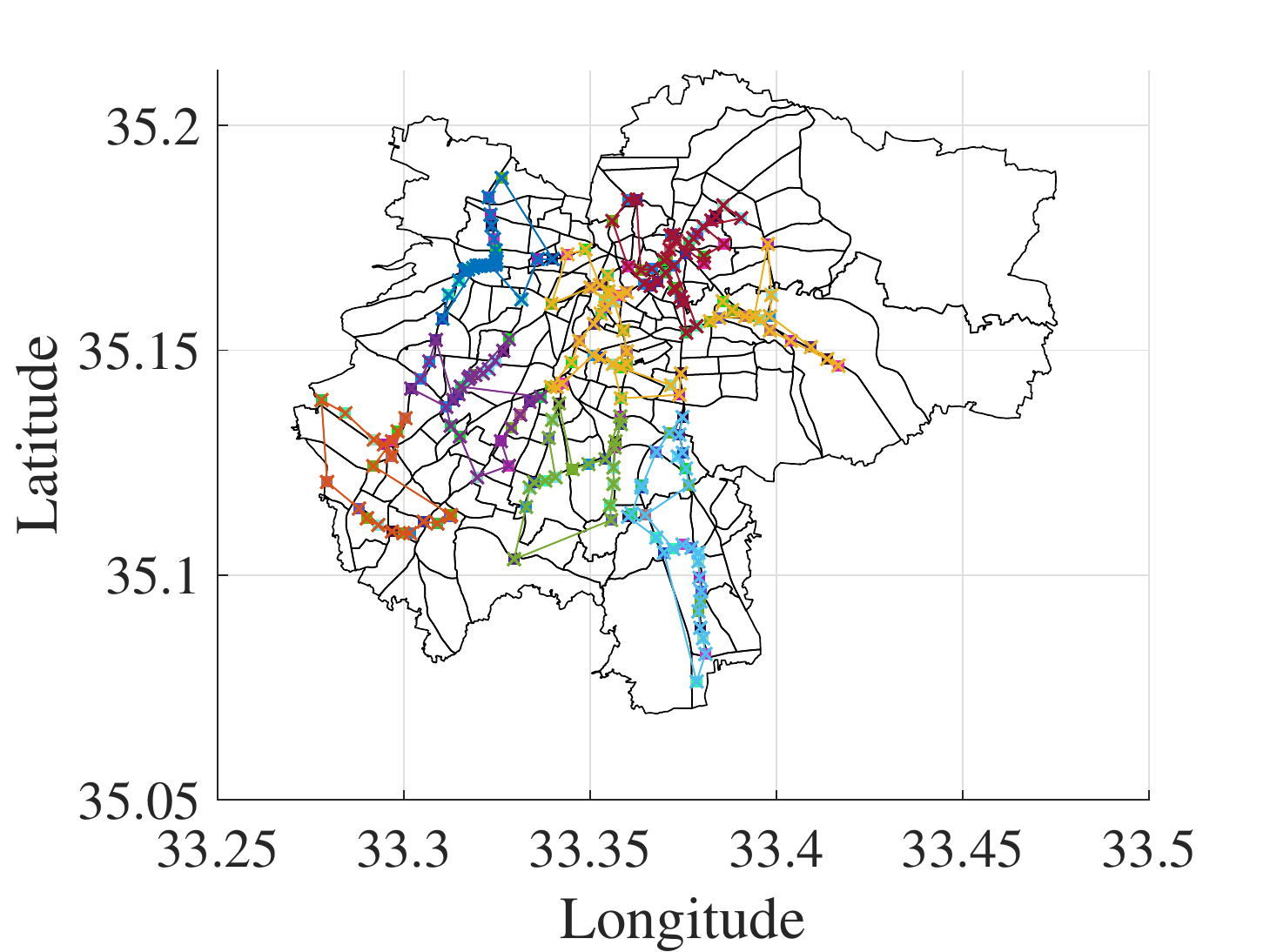}\label{fig:primary_v40_8tours}}
	\subfloat[]{\includegraphics[width=2.4in]{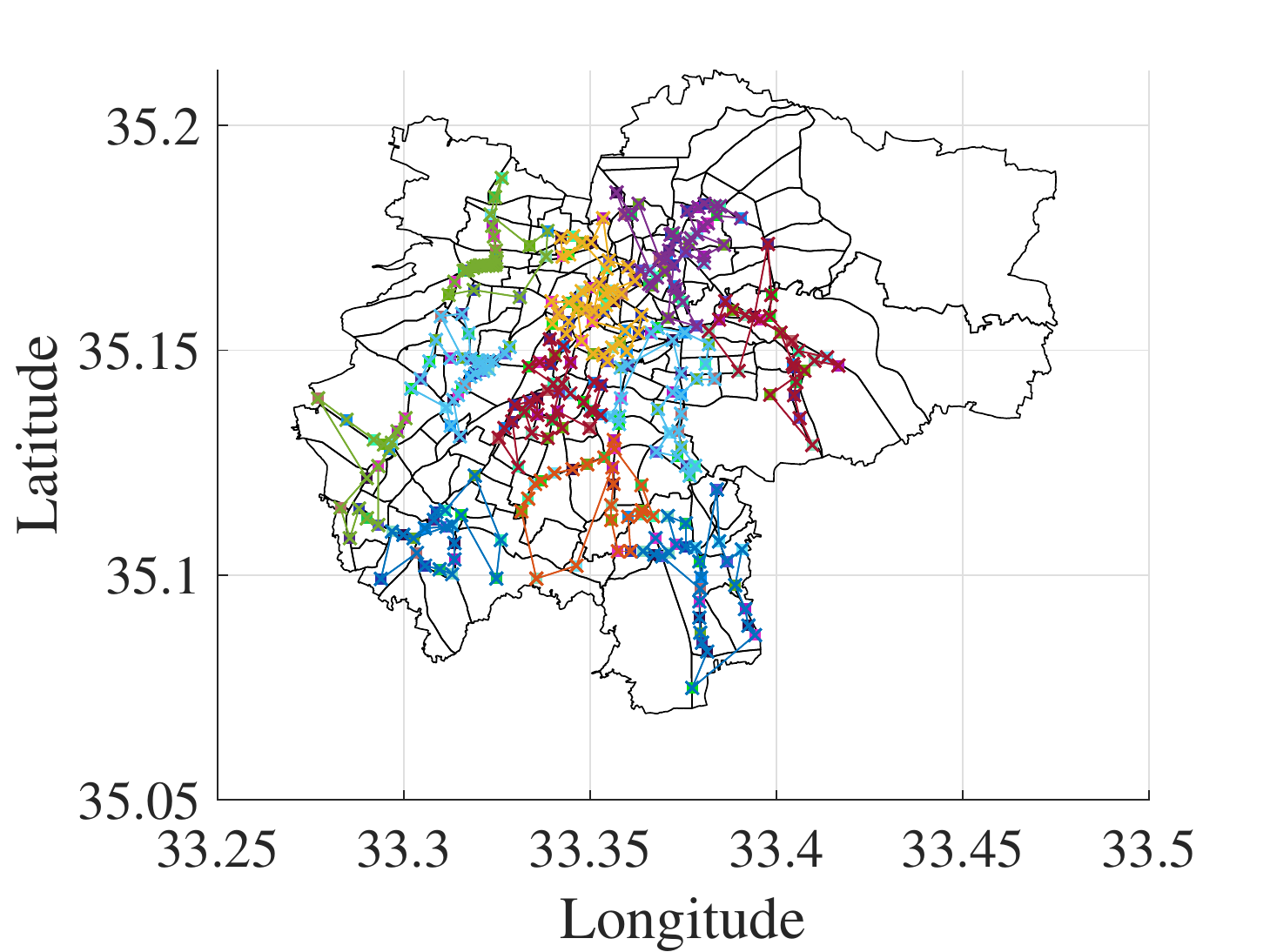}\label{fig:secondary_v40_11tours}}
	\subfloat[]{\includegraphics[width=2.4in]{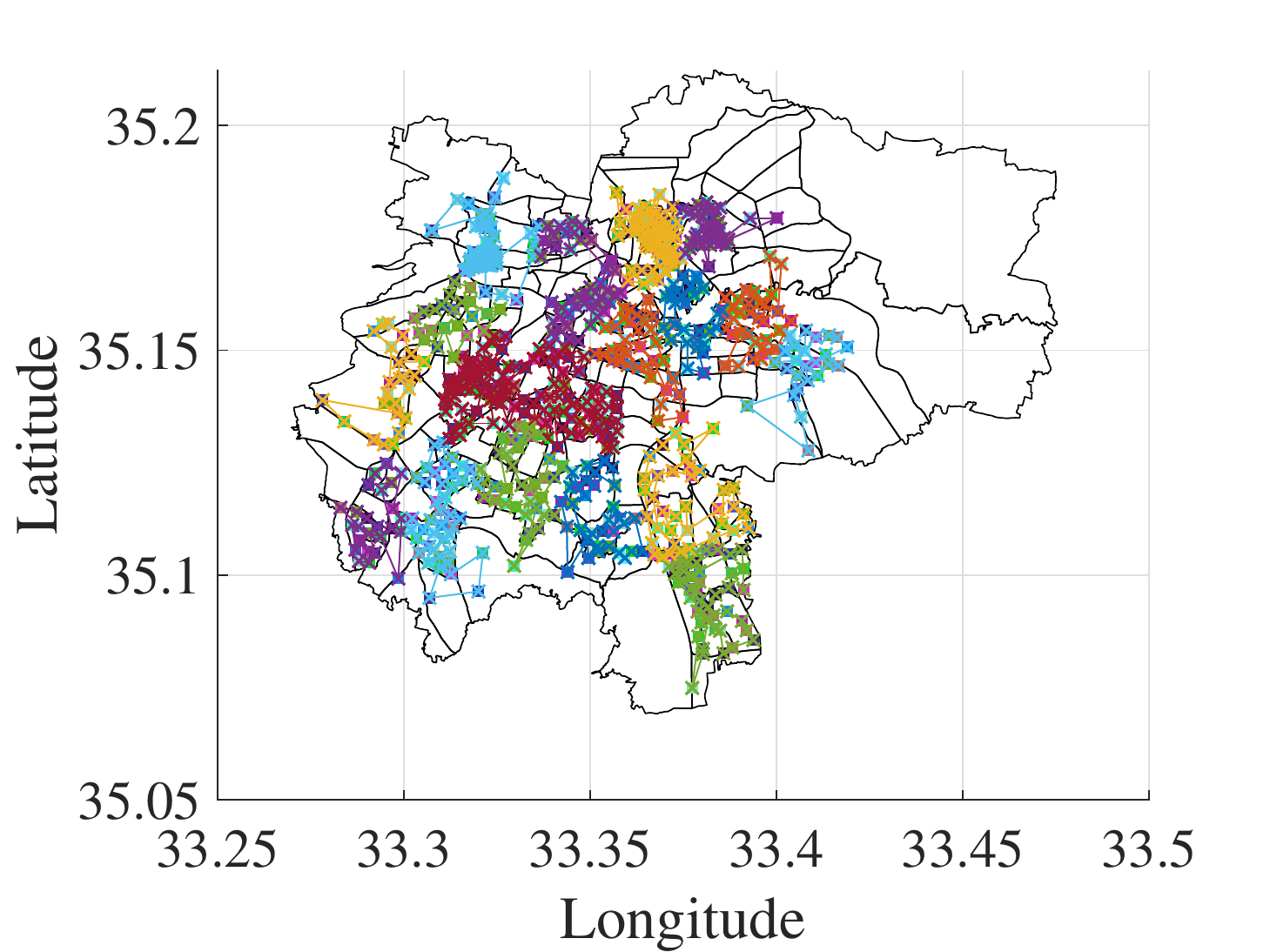}\label{fig:all_v40_18tours}}
	\caption{Tour paths created for (a) the primary road network, (b) primary and secondary network and (c) for all road segments.}
	\label{fig:results}
\end{figure*}

For the underlying road network we have used the following parameters: $h_V=2$ m, $h_U = 500$ m, $w_{BV} = 4$ m. The maximum height of buildings in Nicosia is determined from the postcode (areas of approximately 0.5-2 km$^2$). In the area under consideration, a total of 11 maximum building heights arise: $h_M = \{$5.0,  5.5,  8.0, 10.0, 11.5, 13.5, 17.0, 24.0, 38.0, 45.0, $52.0\}$m which result in the maximum coverage range of $R=\{$664.0, 569.1, 332.0, 249.0, 209.7, 173.2, 132.8,  90.5,  55.3,  46.3,  $39.8 \}$m, respectively. For the monitoring locations we have used a granularity of $\rho = 5$ that provided a total of 86191 candidate locations.

The parameters used for the evaluation of the proposed MTA algorithm are the following: $v= 40$km/h and $t=30$min; thus the maximum travel distance is $m_d=20$ km. 

The proposed framework was tested for the following three network setups:
\begin{itemize}
	\item The primary road network which includes Nicosia’s arterial roads, comprising a set of 198 positions.
	\item The primary and secondary road network which include Nicosia’s main roads in addition to the arterial roads, comprising a set of 305 positions.
	\item The entire road network of Nicosia which includes all roads (primary, secondary and residential), comprising a set of 709 positions. 
\end{itemize}

For the above road networks we have defined a matrix $D$ with the travelling distances between two different positions in the network, by considering a fully connected network. The algorithm was implemented using Matlab.

Fig. \ref{fig:results} depicts the tours developed using MTA for each UAV to monitor the primary, primary and secondary, and the entire road network of Nicosia. In total, 8, 11 and 18 tours were computed, respectively using the aforementioned parameter, demonstrating that with only a very small number of drone units the main metropolitan areas of small-to-medium size cities can be surveyed effectively. 

\begin{figure}
	\centering
	\includegraphics[width=\linewidth]{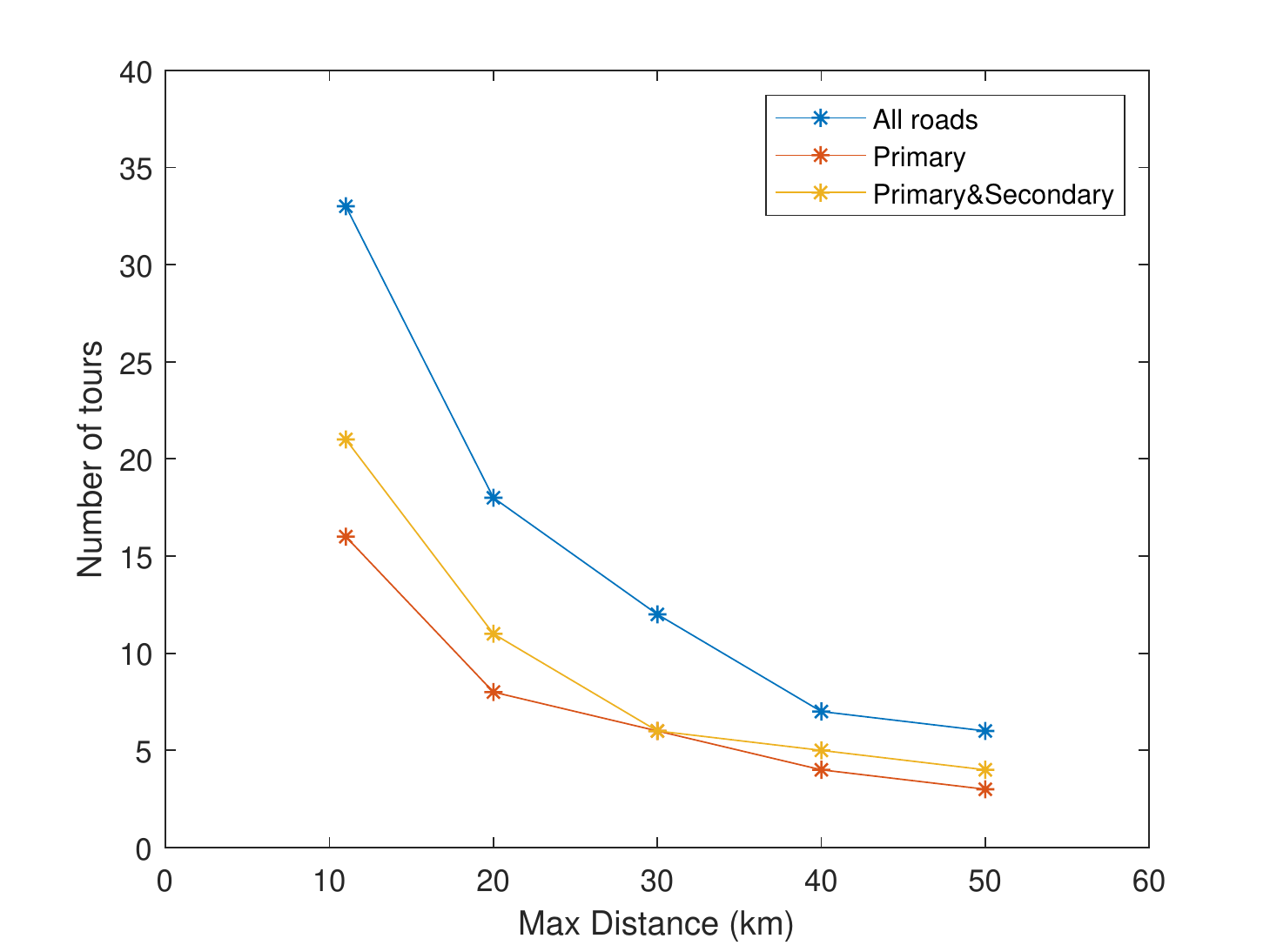}
	\caption{Number of tours created for varying maximum travel distance, $m_d$.}
	\label{fig:velocity_paths}
\end{figure}

To further investigate the impact of flying parameters, we run our algorithm using different velocity values, and thus varying $m_d$. As shown in Fig. \ref{fig:velocity_paths}, there is a non-linear relationship with the number of tours necessary to cover the same area. For instance, take as an example the case of all road segments, for $m_d=40$ instead of $m_d=20$ the number of tours necessary drops from 18 to 7, a $60\%$ reduction in the required tours.

In order to better examine the tour variations between the different clusters, a coefficient of variation is introduced hereafter, as a function of the total travel distance:
\begin{equation}
CV=\frac{std(\{\mathcal{C}_{1},\ldots,\mathcal{C}_{L}\})}{mean(\{\mathcal{C}_{1},\ldots,\mathcal{C}_{L}\})}
\end{equation}
where the standard deviation of the total travelling distance $\mathcal{C}_{l}$ of all clusters $l\in\mathcal{L}$ is divided by the average value. As a standard guideline, a $CV\geq1$ demonstrates a relatively high variation, while a $CV<1$ indicates homogeneity. As shown in table \ref{tab:Var_coef}, the resulting tours experience a low-variance which demonstrates that the computed tours equally contribute to covering the monitoring area and verifies the effectiveness of the proposed algorithm.
\begin{table}[H]
\begin{center}
	\begin{tabular}{|c|c|c|c|}
		\hline
		Velocity (km/h) &	Primary & Primary - Secondary & All Roads \\
		\hline
		20&	0.25& 0.24& 0.23\\
		40& 0.26& 0.29& 0.22\\
		60&	0.14& 0.18&	0.14\\
		80&	0.18& 0.12&	0.18\\
		100& 0.12& 0.10& 0.11\\
		\hline
	\end{tabular}
\vspace{0.1cm}
\caption{CV values between tours created for the three network setups.}
\label{tab:Var_coef}
\end{center}
\end{table}


To demonstrate the applicability of the proposed algorithm, we consider the alternative basic setup where CIA is executed with random initial node placement (RIP). For comparison, the same number of tours is employed for both MTA and RIP. While MTA first clusters nodes and then computes tours, in RIA we simple build tours using random initial nodes for each tour. Table \ref{tab:Var_mean} provides results for 100 Monte Carlo simulations of the RIP variant and its comparison to MTA for the three network setups. The first row indicates the tour length $|T|$ in terms of the number of edges while the second row of the table provides the average tour cost $\bar{\mathcal{C}}_{\mathcal{L}}$ for all computed tours.


\definecolor{cream}{rgb}{1.0, 0.99, 0.82}
\definecolor{azure}{rgb}{0.94, 1.0, 1.0}
\definecolor{palepink}{rgb}{0.98, 0.85, 0.87}
\begin{table}[H]
	\begin{center}
		\newcolumntype{b}{>{\columncolor{azure}}c}
		\newcolumntype{y}{>{\columncolor{cream}}c}
		\newcolumntype{r}{>{\columncolor{palepink}}c}
		\begin{tabular}{|l|r|r|r|y|y|y|b|b|b|}
			\hline
			 &\multicolumn{3}{|r|}{Primary} & \multicolumn{3}{|y|}{Primary-Secondary}& \multicolumn{3}{|b|}{All} \\
			& \multicolumn{1}{|r|}{MTA} & \multicolumn{2}{|r|}{RIP}& \multicolumn{1}{|y|}{MTA} & \multicolumn{2}{|y|}{RIP} & \multicolumn{1}{|b|}{MTA} & \multicolumn{2}{|b|}{RIP} \\	
			& & $\mu$& $\sigma^2$ & & $\mu$ & $\sigma^2$ &  & $\mu$&$\sigma^2$\\
			\hline
			$\!\!\!|T|\!\!\!$&24.8&	24.6&	0.2&	25.4&	27.5&	0.6&	30.8&	39.3&	0.1 \\
			$\!\!\!\bar{\mathcal{C}}_{\mathcal{L}}\!\!\!$&14.2&	17.8&	0.8&	13.8&	17.8&	0.6&	10.7&	17.9&	0.5 \\
			\hline
		\end{tabular}
		\vspace{0.1cm}
		\caption{Comparison of MTA and basic RIP performance for all road network setups.}
		\label{tab:Var_mean}
	\end{center}
\end{table}

Interestingly, the total number of edges included in both the MTA and RIP tours is approximately the same for all three cases. Importantly however, the MTA outperforms the basic RIP approach in all three cases with respect to the total travel cost while also ensuring that the complete road network is monitored. Considering for example the case of monitoring all road segments, the total number of tour legs for MTA is $20\%$ lower than that of RIP and the total cost is close to half of that attained by RIP. These results demonstrate the applicability of the proposed framework and the potential gains in applying this framework in practice.

\section{CONCLUSIONS AND FUTURE WORK}\label{concl}
This work presents a detailed feasibility study in the use of a fleet of drones for road traffic monitoring. To do that, a systematic framework has been derived for quantizing the road network into monitoring locations and using those locations to build tours that a fleet of drones can use to acquire aerial footage of the road traffic conditions. The derived MTA algorithm has been evaluated over a realistic network setup with the road network topology as well as the build-up parameters were extracted from real data. The results demonstrate the applicability of the proposed framework in deriving minimum-cost tours.   

Future work aims at further investigating flying dynamics as well as propulsion energy models in constructing more accurate tour plans as well as investigating how changing altitudes can influence the data acquisition time and accuracy.

\section*{Acknowledgments}
This work is supported by the European Union’s Horizon 2020 research and innovation programme under grant agreement No 739551 (KIOS CoE) and from the Republic of Cyprus through the Directorate General for European Programmes, Coordination and Development.

\end{document}